\documentclass[11pt]{article}
\usepackage{epsfig}
\usepackage{psfrag,graphicx,verbatim,array,multicol,palatino,enumerate}
\usepackage{amsmath,amsfonts,bm,rotating,amssymb,amsthm,latexsym,subfigure}
\usepackage{longtable}
\usepackage{multirow}
\usepackage{authblk}
\usepackage{xcolor}

\makeatletter \@addtoreset{equation}{section} \makeatother

\def\cR{\mathcal{R}}
\def\cL{\mathcal{L}}
\def\cP{\mathcal{P}}
\def\cN{\mathcal{N}}

\def\cX{\mathcal{X}}

\def\cC{\mathcal{C}}
\def\cS{\mathcal{S}}
\def\cech{\v{C}ech }

\def\eps{\epsilon}
\def\iid{\mathrm{i.i.d.}}

\DeclareMathOperator{\PH}{PH}
\DeclareMathOperator{\birth}{birth}
\DeclareMathOperator{\death}{death}
\def\fmin{f_{\min}}

\newcommand{\norm}[1]{\left\|#1\right\|}
\newcommand{\abs}[1] {\left| {#1}\right|}

\newcommand{\pois}[1]{\mathrm{Poisson}\param{{#1}}}
\newcommand{\param}[1]{\left(#1\right)}

\newcommand{\E}{\mathbb{E}} 

\newcommand{\mean}[1] {\E\left\{{#1}\right\}}

\newcommand{\set}[1]{\left\{#1\right\}}

\newcommand{\R}{\mathbb{R}}
\newcommand{\F}{\mathbb{F}}

\newcommand{\prob}{\mathbb{P}}
\newcommand{\expect}{\mathbb{E}}
\newcommand{\var}{\mbox{Var}}

\newcommand{\G}{{G}}
\newcommand{\C}{\mathcal{C}}

\newtheorem{theorem}{Theorem}[subsection]
\newtheorem*{theorem*}{Theorem}
\newtheorem{corollary}[theorem]{Corollary}
\newtheorem{lemma}[theorem]{Lemma}
\newtheorem{proposition}[theorem]{Proposition}

\newtheorem{conj}[theorem]{Conjecture}
\newtheorem*{claim*}{Claim}
\theoremstyle{definition}
\newtheorem{definition}[theorem]{Definition}
\theoremstyle{remark}

\begin{document}

\title{Topology of random geometric complexes: a survey}

\author{Omer Bobrowski}
\affil{Technion - Israel Institute of Technology, Department of Electrical Engineering}

\author{Matthew Kahle}
\affil{The Ohio State University, Department of Mathematics}

%


\maketitle

\baselineskip20pt

\section{Introduction}

In this expository article, we survey the rapidly emerging area of random geometric simplicial complexes. 
Random simplicial complexes may be viewed as higher-dimensional generalizations of random graphs. Perhaps the most studied model of random graph is the Erd\H{o}s--R\'enyi model $\G(n,p)$, where every edge appears independently with probability $p$. Textbooks overviewing this subject include those by Bollob\'as \cite{Bollo} and Janson, {\L}uczak, and Rucinski \cite{randomgraphs}. Simplicial complex analogues of $\G(n,p)$ and their topological properties have been the subject of a lot of activity in recent years. See for example \cite{BHK11,flag,LM,MW} and the references in the survey article \cite{kahle2014topology}.

For certain applications, however, and especially for modeling real-world networks such as social networks, the edge-independent model $\G(n,p)$ is not considered to be particularly realistic. For example, we might expect in a social network that if we know that $X$ is friends with $Y$ and $Z$, then it becomes much more likely than it would be otherwise that $Y$ is friends with $Z$.

Many other models of random graphs have been studied in recent years, and one family of models that has received a lot of attention is the random geometric graphs---see Penrose's monograph \cite{Penrose} for an overview. The random geometric graph $\G(n,r)$ is made by choosing $n$ points independently and identically distributed (i.i.d.), according to a probability measure on Euclidean space $\R^d$ (or any other metric space), and these points correspond to the vertices of the graph. Two vertices $x$ and $y$ are connected by an edge if any only if the distance between $x$ and $y$ satisfies $d(x,y) < r$. Since one is usually interested in asymptotic properties as $n \to \infty$, we usually think of the threshold distance $r$ as a function of $n$.

This is a very general setup, and many variations on this basic model have been studied. The most closely related model to the 
$n$ points i.i.d.\ model is a geometric graph on a Poisson point process with expected number of points $n$. A Poisson point process replaces the independence of points with spatial independence. There is a lot of technology available for transferring theorems between these two models. See, for example, Section 1.7 of \cite{Penrose}.
One might also consider more general point processes than Poisson. For example, Yogeshwaran and Adler \cite{YA12} studied random geometric graphs and complexes over more general stationary point processes. This family includes certain attractive and repulsive point processes, as well as stationary determinantal processes.
In addition, we can consider random geometric graphs in metric measure spaces, such as Riemannian manifolds equipped with probability measures. The topological and geometric properties of such graphs (and their higher-dimensional analogues) were recently studied in \cite{BM14, bobrowski2017random}.

\medskip

There are several natural ways of extending a geometric graph to a simplicial complex, in particular the \cech complex and the Vietoris--Rips complex, whose definitions we review in Section \ref{sec:prelim}.
Our interest in the topology of random geometric complexes will be mainly confined to their homology. Briefly, if $X$ is a topological space, its degree $k$-homology, denoted by $H_k(X)$ is a vector space (assuming field coefficients). The dimension $\dim H_0(X)$ the number of connected components of $X$, and for $k>0$, $H_k(X)$ contains information about $k$-dimensional `holes'. The Betti numbers  of $X$ are defined as $\beta_k(X) = \dim H_k(X)$.

One motivation for studying the topological features of random geometric complexes comes from topological data analysis (TDA). In TDA one builds a simplicial complex (or filtered simplicial complex) on data, and infers qualitative features of the data from homology (or persistent homology) of the point cloud. Studying the topology of random geometric complexes is related to developing probabilistic null hypotheses for topological statistics. We discuss this further in Section \ref{sec:persistent}. The seminal work by Niyogi-Smale-Weinberger \cite{niyogi2008finding, NSW11} introduced a probabilistic analysis to homology recovery algorithms. This was further extended in \cite{balakrishnan_tight_2013, BM14, bobrowski_topological_2014, balakrishnan_statistical_2013}. For surveys of persistent homology in topological data analysis, see Carlsson \cite{Carlsson09} and Ghrist \cite{Ghrist08}.

Studying the limiting behavior of random geometric complexes, the first observation we make is that there exist three main regimes for which the limiting properties of the complexes are significantly different. The term that controls the limiting behavior is $\Lambda = nr^d$, which can be thought of as the average number of points in a ball of radius $r$ (up to a constant).

The subcritical (sometimes called `sparse' or `dust') regime, is when $\Lambda \to 0$. In this regime the geometric complex is highly disconnected, and this is where homology first appears. 

The critical regime (sometimes called `the thermodynamic regime') is when $\Lambda =  \lambda\in(0,\infty)$.  Here, the dimension of homology reaches its peak linear growth, and this is also where percolation occurs (the formation of a `giant' component) --- see the discussion in Section \ref{sec:connectivity_critical}.

Finally, in the super-critical regime we have $\Lambda\to\infty$. In this regime it is known that the number of components slowly decays, until we reach the connectivity threshold. An analogous process occurs for higher homology --- cycles get filled, until eventually every $k$-cycle is a boundary and homology $H_k$ vanishes. But in contrast, for higher homology $k \ge 1$ there is another phase transition where homology $H_k$ first appears.

We note that the connectivity (or $H_0$) properties of random geometric graphs were extensively studied in the past, see \cite{Penrose} for a  comprehensive review. Thus, in this survey we will mainly focus on more recent results related to higher degrees of homology ($H_k,\ k\ge 1$).
 
The rest of this survey is structured as follows. In Section \ref{sec:prelim} we present the concepts and notation that will be used later. Section \ref{sec:connectivity} quickly reviews classical results about the connectivity of random geometric graphs for completeness. Section \ref{sec:betti_numbers} presents a summary of the main results known to date about the limiting behavior of the homology of random geometric complexes. In Section \ref{sec:morse_distance} we review an alternative approach  to study the homology of random \cech complexes using Morse theory for the distance function. Sections \ref{sec:mani} and \ref{sec:stationary} review two extensions to the results in Section \ref{sec:betti_numbers} - one for compact manifolds and the other for stationary point processes. Section \ref{sec:extreme} discusses the case where the distribution underlying the point process has an unbounded support, from an extreme value analysis perspective. In section \ref{sec:persistent} we discuss work in progress that studies the persistent homology generated by random geometric complexes. Finally, in Section \ref{sec:future} we present a list of open problems and future work in this area.

\section{Preliminaries}\label{sec:prelim}

In this section we wish to briefly introduce the concepts and notation that will be used throughout this survey.

\subsection{Homology}\label{sec:homology}
We wish to introduce the concept of homology here in an intuitive rather than a rigorous way. For a comprehensive introduction to homology, see \cite{Hatcher,munkres_elements_1984}.
Let $X$ be a topological space. The \textit{homology} of $X$ is a set of abelian groups $\set{H_k(X)}_{k=0}^\infty$, which are topological invariants of $X$.
 
In this paper we consider homology with coefficients in a field  $\mathbb{F}$, in this case $H_k(X)$ is actually a vector space. The zeroth homology $H_0(X)$ is generated by elements that represent connected components of $X$. For example, if $X$ has three connected components, then $H_0(X) \cong \mathbb{F} \oplus \mathbb{F} \oplus \mathbb{F}$ (here $\cong$ denotes group isomorphism), and each of the three generators corresponds to a different connected component of $X$. For $k\ge 1$, the $k$-th homology $H_k(X)$ is generated by elements representing $k$-dimensional ``holes" or ``cycles" in $X$. An intuitive way to think about a $k$-dimensional hole is as the result of taking the boundary of a $(k+1)$-dimensional body. For example, if $X$ a circle then $H_1(X) \cong \mathbb{F}$, if $X$ is a $2$-dimensional sphere then $H_2(X) \cong \mathbb{F}$, and in general if $X$ is a $n$-dimensional sphere, then
\[
H_k(X) \cong \begin{cases} \mathbb{F} & k=0,n \\
\set{0} & \mbox{otherwise}.
\end{cases}
\]
For another example, consider the $2$-dimensional torus $\mathbb{T}$. The torus has a single connected component  so $H_0({\mathbb{T}}) \cong \mathbb{F}$, and a single $2$-dimensional hole (the void inside the surface) implying that $H_2({\mathbb{T}}) \cong \F$. As for $1$-cycles (or closed loops) the torus has two linearly independent loops, and so $H_1({\mathbb{T}}) \cong \F\oplus \F$.

The dimension of the  $k$-th homology group is called the $k$-th \emph{Betti number},  denoted by $\beta_k(X) := \dim(H_k(X))$.

\subsection{Geometric complexes}

The geometric complexes we will be studying are the \cech and the Vietoris-Rips complexes, defined as follows.
\begin{definition}[\cech complex]\label{def:cech_complex}
Let $\cX = \set{x_1,x_2,\ldots,x_n}$ be a collection of points in $\R^d$, and let $r>0$. The \cech complex $\C_r(\cX)$ is constructed as follows:
\begin{enumerate}
\item The $0$-simplices (vertices) are the points in $\cX$.
\item A $k$-simplex $[x_{i_0},\ldots,x_{i_k}]$ is in $\C_r(\cX)$ if $\bigcap_{j=0}^{k} {B_{r/2}(x_{i_j})} \ne \emptyset$.
\end{enumerate}
\end{definition}

\begin{definition}[Vietoris-Rips complex]\label{def:rips_complex}
Let $\cX = \set{x_1,x_2,\ldots,x_n}$ be a collection of points in $\R^d$, and let $r>0$. The Vietoris--Rips complex $\cR_r(\cX)$ is constructed as follows:
\begin{enumerate}
\item The $0$-simplices (vertices) are the points in $\cX$.
\item A $k$-simplex $[x_{i_0},\ldots,x_{i_k}]$ is in $\cR_r(\cX)$ if $\norm{x_{i_j} - x_{i_l}} \le r$ for all $0\le j,l \le k$.
\end{enumerate}
\end{definition}

\begin{figure}[h!]
\centering
  \includegraphics[scale=0.35]{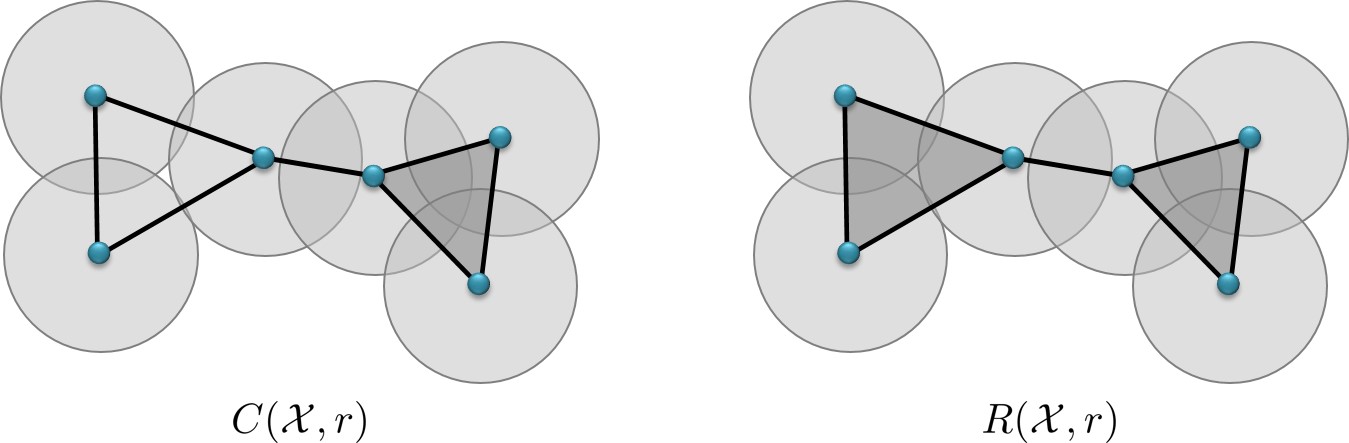}
\caption{On the left - the \cech complex $\C_r(\cX)$, on the right - the Rips complex $R(\cX,r)$ with the same set of vertices and the same radius. We see that the three left-most balls do not have a common intersection and therefore do not generate a 2-dimensional face in the \cech complex. However, since all the pairwise intersections occur, the Rips complex does include the corresponding face.}
\label{fig:cech}
\end{figure}

Figure \ref{fig:cech} shows an example for the \cech and Rips complexes constructed from the same set of points and the same radius $r$, and highlights the difference between them.
As mentioned above, our interest in these complexes will be mostly focused on their homology which is introduced in the next section.

Associated with the \cech complex $\cC_r(\cX)$ is the union of balls used to generate it (in the underlying metric space), which we define as
\begin{equation}\label{eq:union_balls}
	B_{r/2}(\cX) := \bigcup_{x\in \cX}B_{r/2}(x).
\end{equation}
The spaces $\cC_r(\cX)$ and $B_{r/2}(X,r)$ are of a completely different nature. Nevertheless, the following lemma claims that they are very similar in the topological sense.
This lemma is a special case of a more general topological statement originated in \cite{borsuk_imbedding_1948} and commonly referred to as the `Nerve Lemma'.


\begin{lemma}[The Nerve Lemma, Borsuk \cite{borsuk_imbedding_1948}]\label{lem:nerve}
Let $\cC_r(\cX)$ and $B_{r/2}(\cX)$ as defined above. If for every $x_{i_1},\ldots,x_{i_k}$ the intersection $B_{r/2}(x_{i_1})\cap \cdots\cap B_{r/2}(x_{i_k})$ is either empty or contractible (homotopy equivalent to a point), then $\C_r(\cX)\simeq B_{r/2}(\cX)$, and in particular,
\[
	H_k(\cC_r(\cX)) \cong H_k(B_{r/2}(\cX)),\quad \forall k\ge 0.
\]
\end{lemma}
This lemma is highly useful in the study of the random \cech complex, since it allows us to translate questions about the random complex into questions about coverage properties, and enables the use of  Morse theory (see Section \ref{sec:morse_distance}). One immediate consequence of the Nerve Lemma is that if $\cX \subset \R^d$ then $H_k(\C_r(\cX)) = 0$ for all $k\ge d$.

\subsection{Point processes}
Most of the results on random geometric complexes focus on two very similar point processes. In both cases we start with a probability density function $f:\R^d\to \R$, which we always assume to be measurable and bounded.
\begin{itemize}
\item {\bf The binomial process:} \\
 $\cX_n = \{X_1,X_2,\ldots, X_n\}$ is a set of $\iid$ (independent and identically distributed) random variables in $\R^d$ generated by the density function $f$.
 \item  {\bf The Poisson process:}\\
 $\cP_n$ is a spatial Poisson process in $\R^d$ with intensity function $\mu = nf$. The distribution of $\cP_n$ satisfies the following properties:
 \begin{enumerate}
\item For every compact set $A\subset \R^d$ we have $\abs{\cP_n \cap A} \sim \pois{\mu(A)}$, where $\mu(A) = \int_A \mu(x) dx.$	
\item For every two disjoint sets $A,B \subset \R^d$, we have that $\abs{\cP_n\cap A}$ and $\abs{\cP_n \cap B}$ are independent.
\end{enumerate}
 This process is also known as a `Boolean model'.
\end{itemize}
Note $\abs{\cP_n} \sim \pois{n}$, so that $\mean{\abs{\cP_n}} = n$. In addition, given that $\abs{\cP_n} = M$, the process $\cP_n$ consists of $M$ $\iid$ points distributed according to the density function $f$. In other words, the two processes $\cX_n$ and $\cP_n$ are very similar.
We will state most of the results in terms of the binomial process $\cX_n$, and unless otherwise stated, the same results apply to the Poisson process $\cP_n$.

\begin{figure}
\centering
  \includegraphics[width=3in]{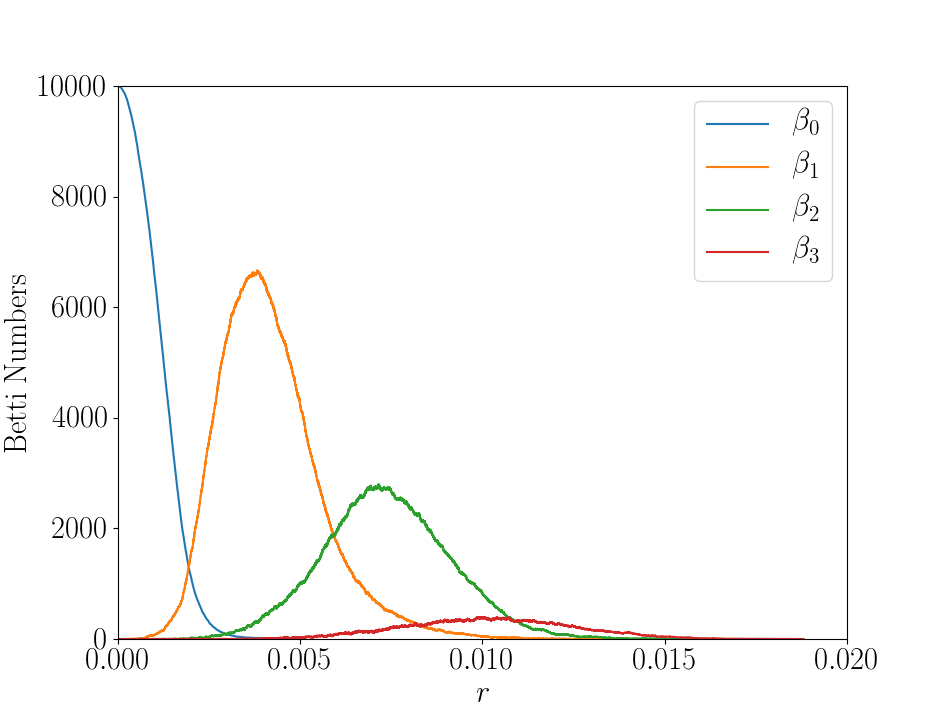}
\caption{The Betti numbers of a random \cech complex as a function of the radius $r$. Here we generated $n=10,\!000$  points uniformly  in $[0,1]^4$. The Betti numbers were calculated using the GUDHI library \cite{gudhi:urm}.}
 \label{fig:rcech}
\end{figure}

In the following we will use the notation $\C_r(n) := \C_r(\cX_n)$, and $\cR_r(n) := \cR_r(\cX_n)$ to state the results about the \cech and Vietoris--Rips complexes generated by the binomial process. Consequently, $\beta_k(n)$ will represent the $k$-th Betti number for either $\C_r(n)$ or $\cR_r(n)$ (which will be clear from the context). Figure \ref{fig:rcech} illustrates the Betti numbers of a random \cech complex, for a fixed $n=10,000$. In most cases we will be interested in the limiting behavior of these complexes as $n\to\infty$ and simultaneously $r = r(n) \to 0$.

\subsection{Convergence of sequences of random variables}

Probability theory uses a number of different notions of convergence. Below we define  the ones used in this survey.

Let $X_1,X_2,\ldots$ be a sequence of real valued random variables, with the cumulative distribution function of $X_n$ given by
\[
	F_n(x) = \prob(X_n \le x),
\]
and let $X$ be a random variable with a cumulative distribution function $F$.
\begin{definition}
$X_n$ converges \emph{in distribution}, or \emph{in law} to $X$, denoted by $X_n \xrightarrow{\cL} X$, if
\[
	\lim_{n\to\infty} F_n(x) = F(x)
\]
for every $x\in \R$ at which $F(x)$ is continuous. 
\end{definition}
This type of convergence is also sometimes referred to as `weak convergence'.

\begin{definition}
$X_n$ converges in $L^p$ to $X$, denoted by $X_n \xrightarrow{L^p} X$, if
\[
	\mean{\abs{X_n-X}^p} \to 0.
\]
\end{definition}
We will mostly use the case $p=2$.

\begin{definition}
$X_n$ converges to $X$ \emph{almost surely}, denoted by $X_n \xrightarrow{a.s.} X$, if
\[
	\prob\left(\lim_{n\to\infty}X_n = X\right) = 1.
\]
\end{definition}

Finally, we have the following probabilistic definition related to limiting events rather than random variables.
\begin{definition}
Let $A_n$ be a sequence of events, perhaps on a sequence of probability spaces. We say that $A_n$ occurs \emph{asymptotically almost surely (a.a.s.)} if
\[
	\lim_{n\to\infty}\prob(A_n) = 1.
\]
\end{definition}

\subsection{Some notation}
Throughout this paper, we use the Landau big-O and related notations. All of these notations are understood as the number of vertices $n \to \infty$. In particular, we write
\begin{itemize}
\item $a_n = O(b_n)$ if there exists a constant $C$ and $n_0 >0$ such that $a_n \le C b_n$ for every $n> n_0$;
\item $a_n = \Omega(b_n)$ if there exists a constant $C$ and $n_0 >0$ such that $a_n \ge C b_n$ for every $n> n_0$;
\item $a_n = \Theta(b_n)$ if  both $a_n = O(b_n)$ and $a_n  = \Omega(b_n)$. We will also denote that by $a_n \sim b_n$;
\item $a_n = o(b_n)$ if $\lim_{n\to\infty} \abs{a_n / b_n} = 0$. We will also denote that by $a_n \ll b_n$;
\item $a_n = \omega(b_n)$ if $\lim_{n\to\infty} \abs{a_n / b_n} = \infty$. We will also denote that by $a_n \gg b_n$.
\end{itemize}
In addition to the above, we use $a_n \approx b_n$ to denote that $\lim_{n\to\infty} a_n/b_n = 1$.

Finally, for any set $A\subset \R^d$ we use $\abs{A}$ to denote the $d$-dimensional volume of the set.

\section{Connectivity}\label{sec:connectivity}
The zeroth homology $H_0$ is generated by the connected components, and its rank $\beta_0$ is the number of components.
Note the connectivity properties of any simplicial complex depend only on its one-dimensional skeleton, namely the underlying graph.
In the \cech and Vietoris--Rips complexes $\C_r(n)$ and $\cR_r(n)$ the underlying graph is the random geometric graph $\G(n,r)$ described above, and therefore the results related to connectivity are the same for both complexes. As we mentioned in the introduction, the main purpose of this survey is to review recent results related to homology in degree $k\ge 1$. However, for completeness, we wish to include a brief review of the key properties related to the connected components. 
Connectivity in graphs is tightly related to the average degree. Note that in the $\G(n,r)$ the degree of a vertex is the number of points lying in a ball of radius $r$ around that vertex. Therefore, for both the binomial and the Poisson processes, the expected degree is proportional to the term
\begin{equation}\label{eq:lambda}
\Lambda := n\cdot r^d.
\end{equation}
As mentioned above, the limiting behavior splits into three main regimes, depending on the limit of the term $\Lambda$. We will correspondingly split the discussion on the limiting results.

\subsection{The subcritical regime}
The subcritical regime (also known as the `sparse' or `dust' regime) is when $\Lambda \to 0$. 
In this regime, the graph $\G(n,r)$ is very sparse, and mostly disconnected. Therefore, the study of connectivity did not draw much attention in the past. See \cite{BM14} for a proof of the following.

\begin{theorem}
If $\Lambda \to 0$ then
\[	
\mean{\beta_0(n)} \approx n.
\]
\end{theorem}

This statement can be sharpened to a central limit theorem, and a law of large numbers can be proved for deviation from the mean. In fact, as we see in the next section, a central limit theorem and law of large numbers continue, even into the critical regime.

\subsection{The critical regime}\label{sec:connectivity_critical}
The critical regime (also known as the `thermodynamic limit') is when $\Lambda = \lambda \in (0,\infty)$. 
In this regime $\beta_0(n) \approx c n $ for some constant $c<1$ (depending on $\lambda$), so the number of components is still $\Theta(n)$, but is significantly lower than in the subcritical regime.
The following law of large numbers is proved in section 13.7 of \cite{Penrose}.

\begin{theorem}[Penrose, \cite{Penrose}]
If $\Lambda = \lambda \in (0,\infty)$, then:
\begin{equation}\label{eq:lim_b0}
	\frac{\beta_0(n)}{n} \xrightarrow{L^2} \int_{\R^d} \param{\sum_{k=1}^\infty k^{-1} p_k(\lambda f(x))} f(x) dx,
\end{equation}
where 
\[
	p_k(t) = \frac{ t^{k-1}}{k!} \int_{(\R^d)^{k-1}}h(0, y_1,\ldots, y_{k-1}) e^{-t A(0,y_1,\ldots,y_{k-1})} dy_1\cdots dy_{k-1},
\]
\[
	h(x_1,x_2,\ldots, x_k) = \begin{cases} 1 & \G(\{x_1,x_2,\ldots, x_k\},1) \textrm{ is connected}, \\ 0 & otherwise,\end{cases}
\]
and
\[
	A(x_1,x_2,\ldots, x_k) := |\bigcup_{j=1}^k B_1(x_j)|.
\]
\end{theorem}
The infinite sum in \eqref{eq:lim_b0} comes from the fact that we need to count the number of components consisting of any possible number of vertices.
The limiting expression provided by the theorem is highly intricate, and at this point impossible to evaluate analytically. 
Nonetheless, as we will discuss later, this theorem provides the only formula available to date for the limit of the Betti numbers in the critical regime.

In addition to a law of large numbers, there is also a central limit theorem available.
\begin{theorem}[Penrose, \cite{Penrose}]
If $\Lambda = \lambda \in (0,\infty)$ then there exists $\sigma>0$ such that 
\[
	\frac{\beta_0(n) - \mean{\beta_0(n)}}{\sqrt{n}} \xrightarrow{\cL} \cN(0,\sigma^2).
\]
\end{theorem}

A more geometric view of connectivity is studied in percolation theory.
Penrose considered the case where $f$ is a uniform probability density on a  $d$-dimensional unit cube, and $\Lambda = \lambda$. A remarkable fact is that there exists a constant $\lambda_c >0$ depending only on the underlying density function, such that if $\lambda < \lambda_c$ then a.a.s.\ every connected component is of order $O( \log n)$, and if $\lambda > \lambda_c$ then a.a.s.\ there is a unique ``giant'' component on $\Theta (n)$ vertices. This sudden change in behavior over a very small shift of parameter is sometimes called a {\it phase transition}.

In chapters 9 and 10 of \cite{Penrose}, Penrose relates percolation on random geometric graphs to more classical {\it continuum percolation theory}. In continuum percolation, also called the Gilbert disk model \cite{Gilbert61}, one considers a random geometric graph on a unit-intensity uniform Poisson process on $\R^d$, and then there is a threshold radius $r_c > 0$ such that for $r > r_c$ the random geometric graph has an infinite connected component, and for $r< r_c$ every component is finite size.
For a deeper study of continuum percolation, see Meester and Roy's book \cite{MR96}. For an introduction and overview of the subject, see Chapter 8 of Bollobas--Riordan \cite{BR06} or Section 12.10 of Grimmett \cite{Grimmett99}.

\subsection{The supercritical regime}\label{sec:connectivity_vanish}
The supercritical regime is when $\Lambda\to \infty$. As we will see soon, if the radius is large enough (yet still satisfying $r\to 0$) then it can be shown that the graph $\G(n,r)$ becomes connected (caveat, this statement depends on the underlying distribution). This phase is sometimes referred to as the `connected regime'. As the radius increases, starting at the critical regime where $\beta_0(n) = \Theta(n)$ and ending at the connected regime where $\beta_0(n) = \Theta(1)$, the number of components in $\G(n,r)$ should exhibit some kind of a decay within the supercritical regime. To this date only partial information is available about this decay process, and we will present it later. We start by describing  the connected regime.

In the case of a uniform distribution on the $d$-dimensional unit box $[0,1]^d$, Penrose gives a sharp result for the connectivity threshold. See \cite{Penrose}, Chapter 13.

\begin{theorem}[Penrose, \cite{Penrose}]  \label{thm:cube}
Let $c \in \R$ be fixed, and set
\[
r = \param{\frac{2^{d-1}}{d\omega_d}\cdot\frac{\log n+c}{n}}^{1/d},
\]
where $\omega_d$ is the volume of the unit ball in $\R^d$.
Then $$\prob( \G(n,r) \mbox{ is connected})\to e^{-e^{-c}}$$ as $n \to \infty$. 
\end{theorem}

\medskip

In other words, the threshold radius for connectivity is $r = \param{\frac{2^{d-1}}{d\omega_d}\cdot\frac{\log n}{n}}^{1/d}$ (or $\Lambda= (2^{d-1}/d\omega_d) \log n$).
It is interesting to contrast Theorem \ref{thm:cube} with the analogous statement for a standard multivariate normal distribution $\cN(0,{\bf{I}}_{d\times d})$ in $\R^d$, a case which Penrose also studies. Here $r$ must be significantly larger, roughly $1 / \sqrt{ \log n}$, in order to ensure connectivity.

\begin{theorem}[Penrose, \cite{Penrose}]\label{thm:gauss}
Let $X_i \sim \cN(0,{\bf{I}}_{d\times d})$ and $c \in \R$ be fixed. If
$$ r=  \frac{(d-1)\log \log n - (1/2) \log \log \log n - 1/\sqrt{4 \pi}  + c}{\sqrt{2 \log n}},$$
then
$$\prob({ G}(n,r) \mbox{ is connected})\to e^{-e^{-c}}$$
as $n \to \infty$.
\end{theorem}
In both cases, letting $c \to \pm \infty$ gives the correct width of the critical window. The critical window is the range of functions $r$ such that the probability of connectedness approaches a constant strictly between $0$ and $1$.

Why does the threshold distance $r=r(n)$ have to be so much larger in the Gaussian case? The support of the Gaussian distribution is unbounded, and there are outlier points at distance roughly $\sqrt{2 \log n}$. The radius must be large enough just to connect these points to the rest of the graph.

The contrast of Theorems \ref{thm:cube} and \ref{thm:gauss} suggests that whatever we hope to prove about the topology of random geometric complexes will necessarily depend on the underlying distribution. On the other hand, certain theorems in geometric probability are fairly general and do not depend on the underlying distribution so drastically.

For example, if we ask what is the threshold for $\G(n,r)$ to contain a given subgraph, or what is the expected number of occurrences  of a given subgraph in the sparse regime, then in some sense the answer does not depend too much on the underlying density function. The following is proved in Chapter 3 of \cite{Penrose}.

\begin{theorem}[Penrose, \cite{Penrose}] \label{thm:subgraph}
Let $\Gamma$ be a finite connected graph on $k$ vertices, and let $N_\Gamma$ count the number of subgraphs isomorphic to $\Gamma$ in $\G(n,r)$. Then 
$$\expect \left[N_\Gamma \right] \sim n^{k} r^{d(k-1)} = n\Lambda^{k-1},$$
as $n \to \infty$. 
\end{theorem}

Note that Theorem \ref{thm:subgraph} applies equally well to uniform distribution on $[0,1]^d$ and to Gaussian distributions; there is no assumption that the underlying measure has compact support. It is only the implied constant in the limit that depends on the measure. This constant may be written out explicitly as an integral - 
\[
	(n\Lambda^{k-1})^{-1} \expect \left[N_\Gamma \right] \approx \frac{1}{k!} \int_{\R^d}f^k(x)dx \int_{\R^{dk}}  h_\Gamma(0, y_1,\ldots, y_{k-1})dy_1\cdots dy_{k-1},
\]
where $h_\Gamma(x_1,\ldots, x_k) = 1$ if $\G(\{x_1,\ldots, x_k\}, 1) \cong \Gamma$ and $0$ otherwise.

As a rule of thumb, one might expect that global properties such as connectivity depend very delicately on the underlying probability measure. Local properties, such as subgraph counts or behavior in the subcritical regime, do not depend so much on the underlying measure.

To conclude this section, we mention a recent result about the supercritical regime preceding connectivity.
As mentioned above, there is a huge gap remaining between the critical regime where $\beta_0(n) = \Theta(n)$ and the connectivity point where $\beta_0(n) = \Theta(1)$. Recent work by Ganesan studies the decay in the number of components within the super critical regime, in the case $d=2$. The assumption is that the underlying probabilty measure on $[0,1]^2$ is supported on a measurable density function $f$, and that $f$ is bounded above and below. The following is Theorem 1 in \cite{Ganesan13}. 

\begin{theorem}[Ganesan, \cite{Ganesan13}]\label{thm:comps}
There exist $a,b,c >0$, such that if $a \log n\le  \Lambda \le b \log n$, then a.a.s.
\[
\beta_0(n) \le n\Lambda^{-1} e^{-c \Lambda},
\]
where the constants $a$ and $b$ depend only on the density function $f$.
\end{theorem}
We will see an analogue of this theorem for higher Betti numbers of the random \cech and  Vietoris--Rips complexes in the following section.

\section{Homology and  Betti Numbers}\label{sec:betti_numbers}
Recall that the $k$-th {\it Betti number} $\beta_k$ is the dimension of $k$-th homology, i.e.\
$$ \beta_k (X) = \dim ( H_k(X) ).$$
As mentioned in the introduction, the homology groups $H_k$ ($k\ge 1$) basically describe cycles (or holes) of different dimensions, and thus the Betti numbers represent the number of cycles.

Betti numbers of random geometric complexes were first studied by Robins in \cite{Robins06}. Robins studies ``alpha shapes'' on random point sets \cite{EKS83}, which are topologically equivalent to \cech complexes but more convenient from the point of view of computation. The underlying distributions are uniform on a $d$-dimensional cube, but to avoid boundary effects periodic boundary conditions are imposed. Robins computes the expected Betti numbers over a large number of  experiments. Furthermore, she explains the shapes of these curves in the ``small radius--low intensity'' regime, writing formulas in the $d=2$ and $d=3$ cases.

The study of the limiting Betti numbers was revisited and significantly extended later in a series of papers by various authors \cite{bobrowski_distance_2011, BM14,geometric, Meckes,  YA12,YSA14}.
In contrast to connectivity which corresponds to reduced zeroth homology $\tilde{H}_0$, the higher homology of random geometric complexes $H_k(\C_r(n))$, $k \ge 1$ is not monotone with respect to $r$. Each homology group passes through two main phase transitions, one where it appears and one where it disappears. 

For the random \cech complex, the phase transition where $H_k$ appears occurs when $\Lambda \sim n^{-\frac{1}{k+1}}$ (or $r \sim n^{-\frac{k+2}{d(k +1)} }$). This radius is within the subcritical regime ($\Lambda\to 0$). In this regime the complex is sparse and highly disconnected which allows very precise Betti number computations --- in  particular we will see that $\beta_k(n) \sim n \Lambda^{k+1}$, and  therefore $\beta_k(n) = o(n)$ 

The phase transition where the $k$-th homology vanishes depends on the underlying probability distribution, but if $f$ has a compact support then we will see that it occurs at $\Lambda = \Theta(\log n)$ (or $r  = \Theta((\log n / n)^{1/d})$), which is within the supercritical regime. This radius is similar to the connectivity threshold we saw in Section \ref{sec:connectivity_vanish}, though the constants are different. The exact vanishing radius for each of the homology groups $H_k$ has not been discovered yet, but it is known that it is controlled by  a second order ($\log\log n$) term that  depends on $k$. We will discuss this  in Section \ref{sec:mani}.

In the critical regime the analysis of the Betti numbers $\beta_k(n)$, $k\ge 1$, is significantly more complicated than the analysis of $\beta_0(n)$. In this case we will see that $\beta_k(n) = \Theta(n)$, however the limiting constants are unknown to date.

We now review the results known to data about the topology of random geometric complexes for each of the regimes.

\subsection{The subcritical regime}\label{sec:betti_subcrit}
The work in \cite{geometric, Meckes} provides a detailed study for the Betti numbers in the subcritical regime. Since a random geometric complex in this regime is so sparse, the vast majority of $k$-cycles are generated by ``small'' sphere-like shapes, with the minimum number of vertices possible. For the \cech complex, the minimum number of vertices to form an $k$-cycle is $k+2$ (for example, to create a $1$-cycle, or a loop, we need at least $3$ vertices). These sphere-like formations are local features, so by the rule of thumb above, we might expect a theorem that holds across a wide class of measures.

A key ingredient in the results is the following indicator function 
\[
	h_k(x_1,\ldots, x_{k+2}) = \begin{cases} 1 & \beta_k(\C_1(\{x_1,\ldots, x_{k+2}\}))=1 \\ 0 & otherwise\end{cases},
\]
testing whether a minimal set forms an $k$-cycle or not. The following theorem provides the limit for the expected Betti numbers.

\begin{theorem} [Kahle, \cite{geometric}]\label{thm:mean_bi_sub}
Let $\Lambda \to 0$, $k\ge 1$ and $d\ge 2$. Then
$$\mean{\beta_k (n) } \approx c_kn \Lambda^{k+1},$$
as $n \to \infty$, where 
\[
	c_k := \frac{1}{(k+2)!} \int_{\R^d}f^{k+2}(x)dx \int_{(\R^d)^{k}}  h_k(0, y_1,\ldots, y_{k+1})dy_1\cdots dy_{k+1}.
\]
\end{theorem}

Theorem \ref{thm:mean_bi_sub} states that $\mean{\beta_k(n)}\sim n\Lambda^{k+1}$. Note that within the subcritical regime the limit of the last term can be either zero, a finite number, or infinity (for different choices of $r$). Combining with the second moment method (see for example Chapter 4 of \cite{Alon}), this is the threshold radius for the phase transition where homology first appears.

\begin{theorem} [Kahle, \cite{geometric}] \label{thm:appear}
Let $d \ge 2$ and $1 \le k \le d-1$ be fixed. Suppose that $\Lambda \to 0$.
\begin{enumerate}
\item If
$$\Lambda  \ll {n^ {-\frac{1}{k +1} }},$$
then a.a.s.\ $H_k( \C_r(n)) = 0$, and 
\item if $$\Lambda \gg {n^{ -\frac{1}{k +1}}}$$
then a.a.s.\ $H_k( \C_r(n)) \neq 0$.
\end{enumerate}
\end{theorem}
Thus, the threshold where the $k$-th homology first appears is $\Lambda = \Theta(n^{-\frac{1}{k+1}})$, or $r = \Theta(n^{ -\frac{k+2}{d(k +1)}})$.

The parallel result for Vietoris--Rips complexes is also given in \cite{geometric}.

\begin{theorem} [Kahle, \cite{geometric}] \label{thm:appear2}
Let $d \ge 2$ and $k \ge 1$ be fixed. Suppose that $\Lambda \to 0$.
\begin{enumerate}
\item If
$$\Lambda  \ll {n^ {-\frac{1}{2k +1} }},$$
then a.a.s.\ $H_k( \cR_r(n)) = 0$, and 
\item if $$\Lambda \gg {n^{ -\frac{1}{2k +1}}}$$
then a.a.s.\ $H_k( \cR_r(n)) \neq 0$.
\end{enumerate}
\end{theorem}

The difference in exponents stems from the fact that in the Vietoris--Rips complex case, the smallest possible vertex support for a nontrivial cycle in $H_k$ is on $2k+2$ vertices (rather than $k+2$ in the \cech complex), a triangulated sphere combinatorially isomorphic to the boundary of the $(k+1)$-dimensional cross polytope.
Another difference is that while in the \cech complex the homology degree is bounded by $d-1$ (a consequence of the Nerve Lemma), for the Vietoris--Rips complex it is unbounded, and we can have cycles of every possible dimension.


Kahle and Meckes studied limiting distributions of Betti numbers in the subcritical regime in \cite{Meckes}.
When $\Lambda = \Theta( n^{-\frac{1}{k+1}})$ (or $r = \Theta(n^{ -\frac{k+2}{d(k+1)}})$), the following is a refinement of Theorem \ref{thm:appear}, and shows that at the threshold where the homology $H_k$ first appears, there is a regime in which the Betti number $\beta_k(n)$ converges in law to a Poisson distribution.

\begin{theorem} [Kahle--Meckes, \cite{Meckes}] \label{thm:kmPoi} Let $1 \le k \le d-1$ and $\mu > 0$ be fixed, and suppose that 
$n\Lambda^{k+1} \to \mu$. Then
$$ \beta_k(n) \xrightarrow{\cL} \pois{\mu c_k},$$
as $n \to \infty$, where $c_k$ is defined in Theorem \ref{thm:mean_bi_sub}.
\end{theorem}

When $r$ is above the threshold, the number of cycles goes to infinity, and with the proper normalization it obeys a central limit theorem. Let $\cN(0,1)$ denote a normal distribution with mean $0$ and variance $1$.

\begin{theorem} [Kahle--Meckes, \cite{Meckes}] \label{thm:kmnor} Let $1 \le k \le d-1$ and suppose that $\Lambda\to 0$
and $$\Lambda \gg {n^ {-\frac{1}{k +1}} }.$$
Then
$$ \frac{ \beta_k(n)  - \expect[ \beta_k(n) ]}{ \sqrt{\var[ \beta_k(n)] }}\xrightarrow{\cL}  \cN(0,1)$$
as $n \to \infty$.
\end{theorem}

Again, because we are in the subcritical regime, these results hold  for a wide variety of measures---whenever the underlying probability measure has a measurable density function which is bounded above. They hold even without compact support, for example for a multivariate normal distribution. In \cite{Meckes} Theorems \ref{thm:kmPoi} and \ref{thm:kmnor} are accompanied by formulas for expectation and variance of the Betti numbers. Parallel limit theorems are also proved for Vietoris--Rips complexes.

\subsection{The critical regime}\label{sec:betti_crit}

The study of the Betti numbers becomes significantly more complicated in the critical regime.
In the subcritical regime, since the random geometric complex is very sparse and disconnected, the vast majority of $k$-cycles are vertex-minimal --- spanning $k+2$ vertices for the \cech complex, $2k+2$ for the Rips. In the critical regime a giant connected component emerges --- see the discussion in Section \ref{sec:connectivity_critical} on percolation theory --- and this significantly complicates the analysis.

To date, there has been some partial progress in studying these cases. For example, we have the following result for expectation.

\begin{theorem} [Kahle, \cite{geometric}] \label{thm:linear}
Suppose that $d \ge 2$ and $0 \le k \le d-1$ are fixed, and $\Lambda = \lambda \in (0,\infty)$. Then for the \cech complex $\cC_r(n)$ we have
$$ \mean{ \beta_k(n)} \sim n.$$
\end{theorem}

A parallel theorem in \cite{geometric} gives the same result for the Vietoris--Rips complex $\cR_r(n)$, but in this case one does not require the assumption that $k \le d-1$; in the critical regime, $\beta_k$ is growing linearly for every $k \ge 0$.

The last theorem provides us with the expected order of magnitude of the Betti numbers, but the actual constants have not yet been discovered. Nevertheless, recent work by Yogeshwaran et al. \cite{YSA14} gives laws of large numbers and central limit theorems for Betti numbers of random \cech complexes in the thermodynamic limit. We state here a few of these results relevant for the \cech complex $\C_r(n)$. The following law of large numbers is Theorem 4.6 in \cite{YSA14}.
\begin{theorem}[Yogeshwaran et al., \cite{YSA14}] 
If $\Lambda = \lambda\in (0,\infty)$, then for each $1\le k \le d-1$ we have almost surely that
\[
\lim_{n\to\infty} \frac{\beta_k(n) - \mean{\beta_k(n)}}{n}  = 0.
\]
\end{theorem}
The version of the central limit theorem proved in \cite{YSA14} is for an underlying uniform distribution, and for simplicity assumes that it is supported on the unit cube in $\R^d$. In this case, they define $I_d(\cP)$ as an interval in $\R$ whose endpoints  are the percolation radii  for $\C_r(n)$ and $\R^d\backslash\C_r(n)$.
\begin{theorem}[Yogeshwaran et al., \cite{YSA14}]
Let $1\le k\le d-1$ and $\Lambda = \lambda\in (0,\infty)$ such that $\lambda \not\in I_d(\cP)$. Then there exists a finite $\sigma^2>0$ such that
\[
\frac{\beta_k(n) - \mean{\beta_k(n)}}{\sqrt{n}} \xrightarrow{\cL} \cN(0,\sigma^2).
\]
\end{theorem}
It is mentioned in \cite{YSA14} that it is not clear whether the restriction to $\lambda \not\in I_d(\cP)$ is required or just a technical artifact of the proof. For the Poisson process $\cP_n$ similar theorems are proved for all $\lambda > 0$.

\subsection{The supercritical regime}\label{sec:betti_supercrit}

In the supercritical regime the correct order of magnitude of the Betti numbers is still not known, but there are bounds. In particular, we have the following for the random Vietoris--Rips complex, which is Theorem 5.1 in \cite{geometric}. 

\begin{theorem} [Kahle, \cite{geometric}]  \label{thm:bound_rips}
Let $\cR_r(n)$ be the random Vietoris--Rips complex, generated by a uniform distribution on a unit-volume convex body in $\R^d$. Then,
\[
\mean{ \beta_k(n)} = O ( n \Lambda^k e^{-c_d \Lambda}),
\]
for some constant $c_d > 0$. Here $c_d$ depends on the dimension $d$ but not on $k$.
\end{theorem}

In particular, if $\Lambda \to \infty$ (the supercritical regime) then $\mean {\beta_k(n)} = o(n)$. 
Theorem \ref{thm:bound_rips} can be compared to Theorem \ref{thm:comps} which bounds the number of connected components.
As an immediate corollary of Theorem \ref{thm:bound_rips} we have the following.

\begin{corollary} \label{cor:van} If $\Lambda \ge c \log n$ then a.a.s. $H_k( \cR_r(n)) = 0$. Here $c$ is any constant such that $c > 1 / c_d$, where $c_d$ is defined in Theorem \ref{thm:bound_rips}. 
\end{corollary}

The proof of Theorem \ref{thm:bound_rips} uses discrete Morse theory to collapse the Vietoris--Rips complex onto a homotopy equivalent CW complex with far fewer faces. 
Combining Theorem \ref{thm:appear2} with Corollary \ref{cor:van} gives the following global picture for vanishing and non-vanishing homology of the random Vietoris--Rips complex.

\begin{theorem} [Kahle, \cite{geometric}] \label{thm:regimes2}
Let $d \ge 2$ be fixed, and suppose that the underlying distribution is uniform on a convex body. Then there exist $a,b$ such that
\begin{enumerate}
\item If
$$\Lambda  \ll n^ {-\frac{1}{2k +1}},$$
then a.a.s.\ $H_k( \cR_r(n)) = 0$, 
\item if $$n^{-\frac{1}{2k +1}}  \ll \Lambda \le a \log n,$$
then a.a.s.\ $H_k( \cR_r(n)) \neq 0$, 
\item and if
$$ \Lambda \ge b \log n$$
then a.a.s.\ $H_k( \cR_r(n)) =0$.
\end{enumerate}
\end{theorem}

For the \cech complex similar bounds are studied in \cite{bobrowski2017random,bob_vanishing}, using Morse theory for the distance function  (discussed in Section \ref{sec:morse_distance}). The idea there is to look for critical points of the distance function, that are responsible for changes in the $k$-th homology.
We note that the following bounds were proven for closed manifolds (compact and without a boundary), while a similar proof can be repeated for the compact and convex case. We shall discuss these bounds in  detail in Section \ref{sec:mani}.

\begin{theorem}  \label{thm:bound_cech}
Let $\C_r(n)$ be the random \cech complex, generated by a uniform distribution on a unit-volume convex body in $\R^d$.  If $\Lambda\to\infty$, then there exist $a_k,b_k >0$ and $c_{d,1},c_{d,2} > 0$ such that
\[
a_k n\Lambda^{k-2}e^{-c_{d,1} \Lambda} \le \mean{ \beta_k(n)} \le b_k n \Lambda^k e^{-c_{d,2} \Lambda}.
\]
\end{theorem}

Combining Theorems \ref{thm:appear} and \ref{thm:bound_cech},   we have the following statement for the \cech complex.

\begin{theorem}[Kahle, \cite{geometric}]  \label{thm:regimes1}
Let $d \ge 2$ and $1 \le k \le d-1$ be fixed, and suppose that the underlying distribution is uniform on a convex body. Then there exist $A,B$ such that 
\begin{enumerate}
\item If
$$\Lambda  \ll n^ {-\frac{1}{k +1}},$$
then a.a.s.\ $H_k( \C_r(n)) = 0$, 
\item if $$ n^ {-\frac{1}{k +1}}  \ll \Lambda \le A \log n,$$
then a.a.s.\ $H_k( \C_r(n)) \neq 0$, 
\item and if
$$ \Lambda \ge B \log n$$
then a.a.s.\ $H_k( \C_r(n)) =0$.
\end{enumerate}
\end{theorem}

Theorems \ref{thm:regimes2} and \ref{thm:regimes1} show that the vanishing threshold radius for higher homology has the same order of magnitude as the  connectivity threshold that we saw in Theorem \ref{thm:cube}, i.e.~it occurs when the average degree is $\Lambda \sim \log n$. Note that this is also when the union of balls $B_{r/2}(\cP_n)$ is known to completely cover the support of the distribution, in which case it can be shown that  $H_k(B_{r/2}(\cP_n)) = 0$. The proof in \cite{geometric} uses this fact together with the Nerve Lemma \ref{lem:nerve} to prove part 3 of the Theorem.  

In Section \ref{sec:mani} we discuss a more refined picture of this transition. We will also see in Section \ref{sec:mani} that these results can be generalized --- for example, to any compact manifold, and for any probability distribution with a density function that is bounded away from zero.

\section{Morse theory for the distance function}\label{sec:morse_distance}
In \cite{bobrowski_distance_2011,bob_vanishing} a different approach was taken to study the homology of \cech complexes which focuses on distance functions.
For a finite set of points $\cP\subset \R^d$ we can define the distance function as follows -
\begin{equation}\label{eq:dist_fn}
	d_{\cP}(x) = \min_{p\in \cP} \norm{x-p}.
\end{equation}
Our  interest in this function stems in the following straightforward observation about the sub-level sets of the distance function:
\[
	d_{\cP}^{-1}([0,\eps]) = B_{\eps}(\cP).
\]
In other words, the sub-level sets of the distance function are exactly the union of balls used to generate a \cech complex. Moreover, from the Nerve Lemma \ref{lem:nerve} we know that these sets have the same homology as the corresponding \cech complex.
Morse theory links the study of critical points of functions with the changes to the homology of their sub-level sets.
Thus, we conclude  that studying the critical points of $d_{\cP}$ might assist us in studying the homology of the \cech complex.
In this section we explore the limiting behavior of the critical points for the random distance function and its consequence to the study of random \cech complexes.

\subsection{Critical points of the distance function}\label{sec:crit_pts}

The classical definition of critical points in calculus is as follows. Let $f:\R^d\to\R$ be a $C^2$ function. A point $c\in \R$ is called a \textit{critical point} of $f$ if $\nabla f (c) =0$, and the real number $f(c)$ is called a \textit{critical value} of $f$. A critical point $c$ is called \textit{non-degenerate} if the Hessian matrix $H_f(c)$ is non-singular. In that case, the \textit{Morse index} of $f$ at $c$, denoted by $\mu(c)$ is the number of negative eigenvalues of $H_f(c)$. A $C^2$ function $f$ is a \textit{Morse function} if all its critical points are non-degenerate, and its critical values are distinct.

Note that the distance function $d_{\cP}$ defined in \eqref{eq:dist_fn} is not everywhere differentiable, therefore the definition above does not apply. However, following \cite{GR97}, one can still define a notion of non-degenerate critical points for the distance function, as well as their Morse index.
Extending Morse theory to functions that are non-smooth has been developed for a variety of applications  \cite{BBK14, Bryzgalova78,GR97, Matov82}. The class of functions studied in these papers have been the minima (or maxima) of a functional
and called `min-type' functions.

\begin{figure}[h]
\centering
 \includegraphics[scale=0.65]{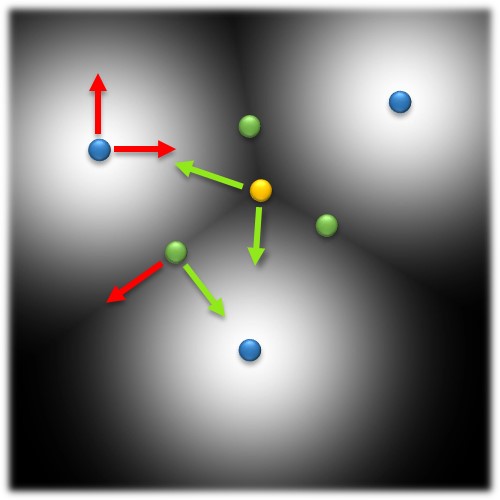}
\caption{\label{fig:crit_pts} Critical points for the distance function in $\R^2$.}
\end{figure}

We wish to avoid the exact definitions of critical points for the distance function and their indexes and introduce them in  a more intuitive way. For the full rigorous definitions and statements see \cite{bobrowski_distance_2011}.
Figure \ref{fig:crit_pts} presents the values of $d_{\cP}$ and the critical points for a set $\cP$ consisting of three points (the blue circles) in $\R^2$.
Obviously, the minima (index $0$ critical points) of $d_{\cP}$ are the points in the set $\cP$ where $d_{\cP} = 0$.
The yellow circle in the middle would be a maximum (index $2$) and the green circles are saddle points (index $1$).
Note that each of the saddle points lie on the segment connecting two sample (blue) points, whereas the maximum lies inside the $2$-simplex spanned by all the three sample points. This is the typical behavior of the critical points of the distance function, and in general we claim that the existence and location of every critical point of index $k$ of $d_{\cP}$ is determined by the configuration of a subset $\cS\subset \cP$ with $\abs{\cS} = k+1$. 

\subsection{Morse Theory}\label{sec:morse}

The study of homology is strongly connected to the study of critical points of real valued functions. The link between them is called Morse theory, and we shall describe it here briefly. For a deeper introduction, we refer the reader to \cite{milnor_morse_1963}.

The main idea of Morse theory is as follows. Suppose that $M$ is a closed manifold (a compact manifold without boundary), and let $f:M\to \R$ be a Morse function.
Denote
\[
M_{\rho} := f^{-1}((-\infty,\rho]) = \set{x\in M : f(x) \le \rho}\subset M
\] (sublevel sets of $f$).
If there are no critical values in $(a,b]$, then $M_a$ and $M_b$ are homotopy equivalent and in particular have isomorphic homology.
Next, suppose that $c$ is a critical point of $f$ with Morse index $i$, and let $v=f(c)$ be the critical value at $c$. Then the homology of $M_{\rho}$ changes at $v$ in the following way. For a small enough $\eps$ we have that the homology of $M_{v+\eps}$ is obtained from the homology of $M_{v-\eps}$ by either adding  a generator to $H_k$ (increasing $\beta_k$ by one) or terminating a generator of $H_{k-1}$ (decreasing $\beta_{k-1}$ by one). In other words, as we pass a critical value, either a new $k$-dimensional cycle is formed, or an existing $(k-1)$-dimensional cycle is bounded or filled.

While classical Morse theory deals with smooth (or $C^2$) Morse functions on compact manifolds \cite{milnor_morse_1963}, it has been extended to many more general situations, and the extenstion to ``min-type'' functions presented in \cite{GR97} enables one to apply similar concepts to the distance function $d_{\cP}$ as well.

Let $\cX_n$ be the binomial process we had before. For $0\le k\le d$, we define  $C_k(r)$ to be the number of critical points of index $k$ of the distance function $d_{\cX_n}$, for which the critical value is less then or equal to $r$. According to Morse theory (and the Nerve Lemma \ref{lem:nerve}), the critical points accounted for by $C_{k}(r)$ are the ones responsible for generating the homology of $\C_r(n)$.

Similarly to the study in Section \ref{sec:betti_numbers} , we can study the limiting behavior of the random values $C_k(r)$ as $n\to\infty$ and $r\to 0$. This was studied in \cite{bobrowski_distance_2011}.
This limiting behavior is in some ways very similar to what we observed for the Betti numbers $\beta_k(n)$. However, as opposed to homology which involves global behavior, the nature of critical points is much more local. This enables us to compute precise limits for $C_k(r)$ even in the critical and super-critical regimes, where the analysis of the Betti numbers at this point has yet to be completed.
We present here the limiting results for the expected values of $C_k(r)$.

\begin{theorem}[Bobrowski--Adler, \cite{bobrowski_distance_2011}]
For $1\le k \le d$ we have,
\begin{enumerate}
\item If $\Lambda \to 0$ then 
\[
	\mean{C_k(r)} \approx \tilde c_k n \Lambda^k;
\]
\item If $\Lambda  = \lambda \in (0,\infty]$ then
\[
	\mean{C_k(r)} \approx \gamma_k(\lambda)\cdot n;
\]
\end{enumerate}
\end{theorem}
The values $\tilde c_k$ and $\gamma_k(\lambda)$ are presented in \cite{bobrowski_distance_2011}, and they depend on the density function $f$, $d$ and $\lambda$ via integration, similarly to the constants $c_k$ in Theorem \ref{thm:mean_bi_sub}.

In the subcritical regime, one can observe that the expected value of $C_k(r)$ is similar to the limit of $\beta_k(n)$ and differs mostly by the index $k$. This is due to the fact that a critical point of index $k$ is generated by a subset of $k+1$ vertices (see discussion above) whereas an $k$-cycle in the subcritical regime is generated by a subset of $k+2$ vertices. Not surprisingly, the distribution of $C_k(r)$ has limit theorems very similar to the ones presented in Section \ref{sec:betti_numbers} for the Betti numbers (see \cite{bobrowski_distance_2011}).

In the critical regime we have $C_k(r) = \Theta(n)$ for all $0\le i \le d$, which, with Morse theory in mind, perfectly agrees with Theorem \ref{thm:linear} stating that $\beta_{k}(n) = \Theta(n)$ as well. As opposed to the Betti numbers, studying the critical points yields precise limits for the expectation as well as a central limit theorem (cf. \cite{bobrowski_distance_2011}). This will enable us later to get a very interesting conclusion regarding the Euler characteristic of $\C_r(n)$.

In the super-critical regime, we still have the exact limits for the number of critical points. However, in this case, it will not reveal much information about $\C_r(n)$, since most of the critical points accounted for by $C_k(r)$ were formed in the critical regime (note that $C_k(r)$ is a monotone function of $r$), and the number of critical points actually being formed in the super-critical regime is actually $o(n)$.
Nevertheless, in some cases (see Section \ref{sec:mani}), it is possible to study the behavior of critical points within the super-critical regime in a finer resolution and use that to draw conclusions about the vanishing of the different degrees of homology.

\subsection{The Euler characteristic}

The Euler characteristic of a simplicial complex $\cS$ has a number
of equivalent definitions, and a number of important applications. One of
the definitions,  via Betti numbers, is
\begin{equation}\label{eq:EC_betti}
    \chi(\cS) = \sum_{k=0}^{\infty} (-1)^k\beta_k(\cS).
\end{equation}
Thus, one can think of the Euler characteristic as an integer ``summary'' of the set of Betti numbers of the complex. 
In the case of the random \cech complex $\C_r(n)$ we  have
\[
	\chi_r(n) := \chi(\C_r(n)) =  \sum_{k=0}^d (-1)^k \beta_{k}(n).
\]
However,  using Morse theory for the distance function, $\chi_r(n)$ can also be computed in the following way
\[
	\chi_r(n) := \sum_{k=0}^d (-1)^k C_k(r).
\]

The limiting behavior of the critical points presented in Section \ref{sec:morse}, thus leads us to the following conclusion.

\begin{corollary}[Bobrowski--Adler, \cite{bobrowski_distance_2011}]\label{cor:cech_ec}
Let $\chi_r(n)$ be the Euler characteristic of $\C_r(n)$, and let $\Lambda = \lambda\in(0,\infty)$. Then
\begin{equation}\label{eq:cech_ec}
\lim_{n\to\infty} n^{-1} \mean{\chi_r(n)}=1+\sum_{k=1}^d {(-1)^k\gamma_k(\lambda) },
\end{equation}
where $\gamma_k(\lambda)$ are increasing functions of $\lambda$ and are defined in \cite{bobrowski_distance_2011}.
\end{corollary}

Note that \eqref{eq:cech_ec} cannot  be proven using only the
existing results on Betti numbers, since  the values of the limiting mean in
the critical regime are not available.
This demonstrates one of the advantages of studying the homology of the
\cech complex via the distance function. An alternative way to compute the Euler characteristic is
\[
\chi_r(n) = \sum_{k=0}^\infty (-1)^k \Delta_k(r),
\]
where $\Delta_k(r)$ is the number of $k$-simplexes in $\C_r(n)$. In \cite{ferraz2011statistics} the Euler characteristic was studied this way for a uniform distribution on a $d$-dimensional torus. Computing the mean value (and also the variance) of $\Delta_k(r)$ is possible, however there are going to be infinitely many summands in this formula, which will make the it highly complicated. Thus, counting critical points is still advantageous.

\begin{figure}[h!]
\centering
 \includegraphics[trim=27cm 6cm 0cm 5cm, clip=true, width=7cm]{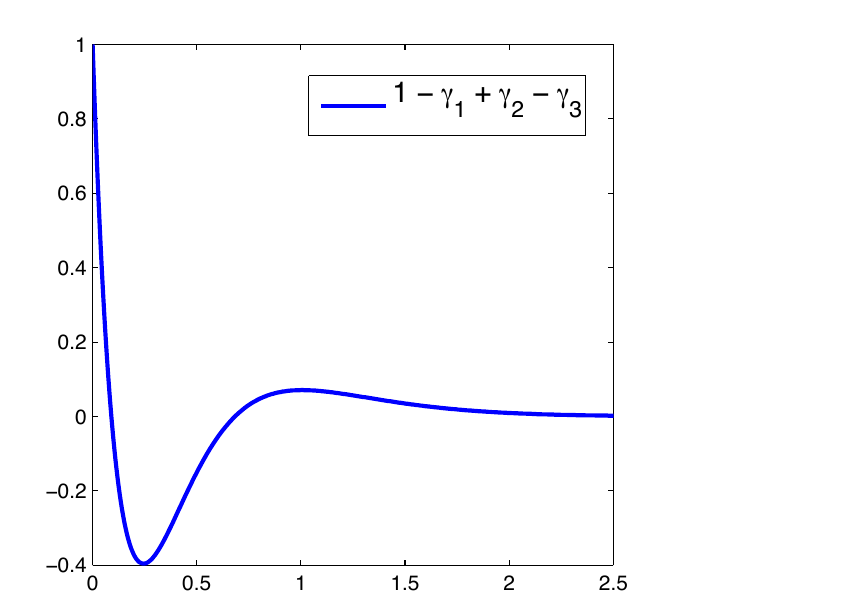}
\caption{\label{fig:ec_crit} The limiting Euler characteristic curve for a uniform distribution on the unit cube in  $\R^3$.}
\end{figure}

Figure \ref{fig:ec_crit} presents the limiting expected Euler characteristic (divided by $n$) as a function of $\lambda$ for a uniform distribution on the unit cube in $\R^3$. In this case the functions $\gamma_k$ ($k=1,2,3$) were computed explicitly in \cite{BM14}) and are given by - 
\begin{align*}
	\gamma_1(\lambda) &= 4(1-e^{-\frac{4}{3}\pi \lambda}),\\
	\gamma_2(\lambda) &= (1+\frac{\pi^2}{16}) (3- 3e^{-\frac{4}{3}\pi \lambda}-4\pi\lambda e^{-\frac{4}{3}\pi \lambda}),\\
	\gamma_3(\lambda) &=\frac{\pi^2}{48} (9 - 9e^{-\frac{4}{3}\pi \lambda} - 12 \pi \lambda e^{-\frac{4}{3}\pi \lambda} - 8\pi^2 \lambda^2e^{-\frac{4}{3}\pi \lambda}).
\end{align*}
Note that the curve starts at positive values, turns negative and then becomes positive once and for all. In $\R^3$ the formula \eqref{eq:EC_betti} implies that $\cX = \beta_0-\beta_1+\beta_2$.

The shapes of the Betti number curves in Figure \ref{fig:rcech} suggests the conjecture that each of the different Betti numbers becomes dominant in a slightly different regime. A similar phenomenon is known to occur for certain random abstract simplicial complexes \cite{flag}, but it is still not known whether this holds for random geometric complexes for the Rips complex as well.

\section{Extending to manifolds} \label{sec:mani}

In sections \ref{sec:connectivity}-\ref{sec:morse_distance} the distributions studied are supported on $d$-dimensional subsets of $\R^d$. The work in \cite{BM14} studied the same type of problems for the case where the distributions are supported on a closed $m$-dimensional manifold embedded in $\R^d$ ($m<d$). In \cite{bob_vanishing} the flat torus was studied as a special case of a Riemannian manifold, and this was extended later to compact (smooth) Riemannian manifolds in \cite{bobrowski2017random}. In this section we will limit the discussion to the \cech complex, although some of the results (in particular the behavior in the subcritical and critical regimes) could be similarly generalized.

\subsection{Closed Manifolds Embedded in $\R^d$}

The exact setup studied was as follows. Let $M\subset \R^d$ be a $m$-dimensional smooth closed manifold (compact and without a boundary). Let $f:M\to \R^d$ be a probability density function supported on $M$. Let $\cX_n = \{X_1,\ldots, X_n\}$ be a set of $\iid$ points generated by $f$, and let $\C_r(n)$ be the \cech complex generated by these points (using $d$-dimensional balls). The results in this case turn out to be very similar to the ones we described earlier, even though the proofs require different analysis tools. In the following we briefly review the results in \cite{BM14} and highlight the main difference from the results in $\R^d$.

The first thing to note is that here, the average degree behaves like $\Lambda = nr^m$ ($m$ being the intrinsic dimension of the manifold). In the subcritical regime, the results for both the Betti numbers $\beta_k(n)$ and the number of critical points $C_k(r)$ are almost identical to those presented in Sections \ref{sec:betti_subcrit} and \ref{sec:morse}. The main difference is that the ambient dimension $d$ is replaced by the intrinsic dimension $m$, and the limiting constants are a bit different. For example, we have that 
$$\mean{\beta_k(n)} \approx \tilde c_k n \Lambda^{k+1},$$
where
$$	\tilde c_k = \frac{1}{(k+2)!} \int_{M}f^{k+2}(x)dx \int_{(\R^m)^{k}}  h_k(0, y_1,\ldots, y_{k+1})dy_1\cdots dy_{k+1}.$$
 These differences stem from the fact that in the subcritical regime the Betti numbers computation is very `local', and locally, a $m$-dimensional manifold looks very similar to $\R^m$. In the critical regime we also have very similar statements to the Euclidean setup.

The main difference in studying manifolds shows up when we study the vanishing of the homology. When studying compact and convex bodies, Theorem \ref{thm:regimes1} states that homology completely vanishes when $\Lambda \sim \log n$ (or $r \sim \left(\frac{\log n}{n}\right )^{1/d}$). Sampling from a manifold, by the Nerve Lemma, we expect that upon coverage the homology of the complex $\C_r(n)$ will not  vanish but rather become equal to the homology of $M$. This result is stated in the following theorem.

\begin{theorem}[Bobrowski--Mukherjee, \cite{BM14}] \label{thm:lim_betti_sup}
Let $\epsilon > 0$ be fixed.
If 
\[
\Lambda \ge \param{\frac{2^m}{\omega_m \fmin}+\eps} \log n
\]
then $H_k(\C_r(n))\cong H_k(M)$ for all $0\le k \le m$ a.a.s., and if
\[
\Lambda \le \param{\frac{2^m}{\omega_m \fmin}-\eps} \log n
\]
then $H_k(\C_r(n))\not\cong H_k(M)$ for all $1\le k \le m$ a.a.s., where $\omega_m$ is the volume of the $m$-dimensional unit ball, and
$f_{\min} = \inf_{x\in M} f(x) >0$.
\end{theorem}
We note that while the second part of this theorem did not appear explicitly in \cite{BM14}, it is a direct consequence of the calculations done there in addition to the Morse theoretical arguments made in \cite{bob_vanishing} (discussed later).
Also note that the vanishing radius for $H_k$ ($k\ge 1$) is twice the radius of connectivity in the same setup (an analog result of Theorem \ref{thm:cube} was proved for the flat torus in \cite{Penrose}, and can be extended to any compact embedded or Riemannian manifold using the techniques in \cite{BM14, bobrowski2017random}).
This phenomenon has a non formal, yet convincing, explanation. In \cite{Penrose} (Theorem 13.17) it is shown that at the edge of connectivity the graph $\G(n,r)$ consists roughly of a giant component and some isolated vertices. For a vertex to be isolated, a ball of radius $r$ around it has to be vacant (i.e. with no other points in $\cX_n$ inside it). To get all the higher homology groups correctly, we need to guarantee that the balls of radius $r/2$ (the ones used to construct the \cech complex) cover the support. Now, the support is covered if and only if there is no vacant ball of radius $r/2$. Thus, it seems harder to reach coverage than connectivity, and the vacancy radii involved have the same ratio as the thresholds we presented.

The statement in Theorem \ref{thm:lim_betti_sup} has an important consequence to problems in manifold learning, since it shows that by studying \cech complexes we can recover the homology of an unknown manifold $M$ from a finite (yet probably large) number of random samples. The analysis of this type of ``topological manifold learning" was established by the seminal work in \cite{niyogi2008finding} and \cite{NSW11}, and Theorem \ref{thm:lim_betti_sup} can be viewed as an  asymptotic and extended version of the main results there. Considering asymptotic behavior has the advantage of covering a more general class of distributions and using fewer assumptions.

Theorem \ref{thm:lim_betti_sup} shows that for large enough radii, the Betti numbers computed $\beta_k(n)$ converge to the Betti numbers of the manifold $\beta_k(M)$. Denoting the error by
$$\hat\beta_k(n)  = \beta_k(n) - \beta_k(M),$$
Theorem \ref{thm:lim_betti_sup} can be viewed as describing the vanishing of the `noisy homology' (so that  $\hat \beta_k(n) \to 0$).

\subsection{Riemannian Manifolds and Homological Connectivity}

The work in \cite{bobrowski2017random,bob_vanishing} studied a similar case to the previous one, only that now the random point process is generated on a $d$-dimensional Riemmanian manifold $(M,g)$. The main difference in this setup, is that now the balls used to create the geometric complexes, are $d$-dimensional intrinsic balls on the manifold (i.e.~using the Riemannian rather than the Euclidean metric). As before, most of the statements we had for random geometric complexes in Euclidean spaces, can be extended to the Riemannian setting. In this section we focus on one particular aspect that has been further studied in the case of compact Riemannian manifolds. In the following we will limit ourself to uniform distributions on manifolds with a unit volume (in which case $f \equiv 1$).

By `homological connectivity' we refer to the phenomenon described above where the $k$-th homology of the \cech complex becomes isomorphic to that of the underlying manifold (i.e.~$H_k(\C_r(n))\cong H_k(M)$). We note that this term was coined by Linial and Meshulam in \cite{LM}.  The result in Theorem \ref{thm:lim_betti_sup} (which could be extended to compact Reiamannian manfiolds) states that for all $k\ge 1$ homological connectivity for $H_k$ occurs around $\Lambda = (2^d/\omega_d) \log n$.
Note, however, that this result does not differentiate between the different homology groups.
Since our previous study shows that cycles in different dimensions are formed by different type of structures, and occur at different radii, we also expect to observe differences in the homological connectivity thresholds for different dimensions $k$.

The work in \cite{bob_vanishing} revisited the study of critical points for the distance function for the case when $M$ is the flat torus (i.e.~$\mathbb{T}^d = [0,1]^d \backslash \set{0\sim 1}$). By providing more details estimates to the number of critical points, the following statement was proved.

\begin{proposition}[Bobrowski \& Weinberger, \cite{bob_vanishing}]  \label{thm:bound_cech_m}
Let $1\le k \le d-1$. If $\Lambda\to\infty$, then there exist $a_k,b_k >0$ such that
\[
a_k n\Lambda^{k-2}e^{-2^{-d}\omega_d \Lambda} \le \mean{ \beta_k(n)} \le \beta_k(M)+ b_k n \Lambda^k e^{-2^{-d}\omega_d \Lambda}.
\]
\end{proposition}
To get the upper bound, we denote by $\hat C_k(r)$ the number of critical points whose critical value is bigger than $r$.
Then $\beta_k(n) \le  \beta_k(M) +\hat C_{k+1}(r)$ since by Morse theory all the cycles in $H_k(\C_r(n))$ that do not belong to $H_k(M)$ are to be  terminated by some critical point of index $k+1$. For the lower bound, we look for critical points of index $k$ with a special local behavior that guarantees to generate a new $k$-cycle (See \cite{bob_vanishing} for details). The last inequality then leads to the following result.

\begin{theorem}
Let $1\le k \le d-1$ and suppose that $w(n)\to \infty$ as $n\to\infty$. Then,
\[
\lim_{n\to\infty} \prob({H_k(\C_r(n))\cong H_k(\mathbb{T})}) = \begin{cases} 1 & \Lambda = (2^d/ \omega_d)(\log n + k\log\log n + w(n)),\vspace{5pt}\\
0 & \Lambda = (2^d/\omega_d)(\log n + (k-2)\log\log n - w(n)).\end{cases}
\]
\end{theorem}
Note that: (a) This statement is about isomorphism of the homology groups, which is stronger than just the equality of the Betti numbers; (b) There is a gap in this description of the phase transition, as the two thresholds differ a $\log\log n$ factor.
In \cite{bobrowski2017random} these results were extended from the flat torus to any compact smooth $d$-dimensional Riemannian manifold. However,  it is not clear how this result generalizes to spaces that have boundaries (as the ones in Section \ref{sec:betti_supercrit}).

Finally, we note that we believe the following conjecture to be the most accurate description of the phase transition for homological connectivity.
\begin{conj}
Let $(M,g)$ be a smooth $d$-dimensional compact Riemannian manifold.
Let $1\le k \le d-1$ and suppose that $w(n)\to \infty$ as $n\to\infty$. Then,
\[
\lim_{n\to\infty} \prob({H_k(\C_r(n))\cong H_k(M)}) = \begin{cases} 1 & \Lambda = (2^d/\omega_d)(\log n + (k-1)\log\log n + w(n)),\vspace{5pt}\\
0 & \Lambda = (2^d/\omega_d)(\log n + (k-1)\log\log n - w(n)).\end{cases}
\]
\end{conj}
The reason why this conjecture should be true is that the same phase transition can be shown to describe the vanishing of  isolated $k$-faces ($k$-simplexes that do not have any $(k+1)$-coface). In all other random simplicial complexes studied in he past it was shown that these isolated faces generate the last cycles that prevent homology from converging. Proving this conjectures, however, remains as future work.

\section{Stationary point processes}
\label{sec:stationary}
The results we presented so far in this survey describe the behavior of geometric complexes constructed from either the binomial process $\cX_n$ or the Poisson process $\cP_n$. Both models exhibit a strong level of independence which plays a significant role in the proofs. For the binomial process $\cX_n$ the number of points is fixed, while the locations of the points are independent. For the Poisson process $\cP_n$ the amount of points in different regions are independent, and given the number of points in a region their locations are independent. 

Recent work by Yogeshwaran and Adler \cite{YA12} extends some of the results presented in this survey to a more general class of spatial point processes allowing certain attractive and repulsive point processes, as well as stationary determinantal processes. In this section we wish to briefly review their results.

A general point process in $\R^d$ can be thought of as a random measure $\Phi(\cdot) = \sum_i \delta_{X_i}(\cdot)$ where $\delta_x$ is the Dirac delta measure concentrated at $x$. In that case, for every subset $A\subset\R^d$, $\Phi(A)$ is a random variable counting the number of points lying inside $A$. The distribution of a random point process $\Phi$ can be characterized by its factorial moment measure functions $\alpha^{(m)}$ defined as follows - 
\[
\alpha^{(m)}(B_1, \ldots, B_m) = \mean{\prod_{i=1}^m \Phi(B_i)},
\]
where $B_1,\ldots, B_m$ are disjoint Borel subsets of $\R^d$. A \emph{stationary} point process is such that the functions $\alpha^{(m)}$ are translation invariant. For example, for the homogeneous Poisson process with constant rate $\mu$, we have that
\[
	\alpha^{(m)}(B_1,\ldots, B_m) = \lambda^m \prod_{i=1}^k \abs{B_i},
\]
which depends only on the volumes of the sets and therefore invariant to translations.
Note that if $\Phi$ is a stationary point process, and $\C_r(\Phi)$ is the corresponding \cech complex, then depending on $r$ either $\mean{\beta_k(\C_r(\Phi))} = 0$ or $\mean{\beta_k(\C_r(\Phi))} = \infty$ (since the process is supported in an infinite domain). Therefore, it does not make sense to try to analyze $\beta_k(\C_r(\Phi))$. Instead, we can define
\[
	\Phi_n := \Phi \cap \left[\frac{-n^{1/d}}{2}, \frac{n^{1/d}}{2}\right],
\]
and try to study
\[
	\beta_k^\Phi(n) := \beta_k(\C_r(\Phi_n)).
\]
Note that if $\Phi$ is a homogeneous Poisson process with rate $\mu =1$, and $\cP_n$ is the Poisson process we used previously supported on the unit cube, then $\C_r(\Phi_n)$ is a scaled version of $\C_{n^{-1/d}r}(\cP_n)$, and so $\beta_k^\Phi(n) = \beta_k(n)$. Therefore, we can view the results in \cite{YA12} as an extension of the models described earlier in this survey. Similarly to the study of the binomial and the Poisson processes we described before, the limiting behavior of $\beta_k^\Phi(n)$ splits into three main regimes. Due to the different scaling, the term controlling the limiting behavior is $r$ rather than $\Lambda$. 

The sparse (or the subcritical) regime is when $r\to 0$. In this case, \cite{YA12} shows that there exists a sequence of functions $f^k$ such that either $f^k\equiv 1$ or $\lim_{r\to 0}f^k(r) = 0$ (depending on the distribution of $\Phi$), and then $$\mean{\beta_k^\Phi(n)} \sim nr^{d(k+1)}f^{k+2}(r),$$
where the exact limiting constant is given by a formula similar in spirit to $c_k$ in Theorem \ref{thm:mean_bi_sub}. The results in \cite{YA12} also provide equivalent limits for the distribution as in Theorems \ref{thm:kmPoi}-\ref{thm:kmnor}.

The critical (thermodynamic) regime is when $r = \lambda \in (0,\infty)$. In this case, \cite{YA12} shows that $\mean{\beta_k^\Phi(n)} = \Theta(n)$ and provide a limit for the Euler characteristic similarly to Corollary \ref{cor:cech_ec}.

Finally, in the super critical regime ($r\to \infty$) \cite{YA12} discusses the connectivity regime, which is when $r^d = \Theta(\log n)$.
Similarly to Theorem \ref{thm:regimes1} they show that there exists a constant $c$ such that if $r \ge c\param{\frac{1}{\log n}}^{1/d}$ then $\C_r(\Phi_n)$ is a.a.s contractible.

In addition to the Betti numbers of the \cech complex, they also provide equivalent results for the Vietoris-Rips complexes $\cR_r(\Phi_n)$ and for the critical point counts $c_k$ for the distance function $d_{\Phi_n}$. In \cite{YSA14} these theorems are extended in some cases, to laws of large numbers and central limit theorems.

\section{Extreme value analysis of random geometric complexes}\label{sec:extreme}
The results in the supercritical regime ($\Lambda\to\infty$) that we presented so far, assumed that the point process is  generated by a distribution with a bounded support (see e.g.~Theorems \ref{thm:cube}, \ref{thm:bound_rips},\ref{thm:bound_cech}). 
As the result in Theorem \ref{thm:gauss} suggests, the limiting behavior can be significantly different once we generate the point process by a distribution with an unbounded support (e.g.~the Gaussian distribution). The work in \cite{adler2013crackle,owada_limit_2015} studied the distribution of the Betti numbers in these cases.

\begin{figure}
\begin{center}
\includegraphics[scale=0.4]{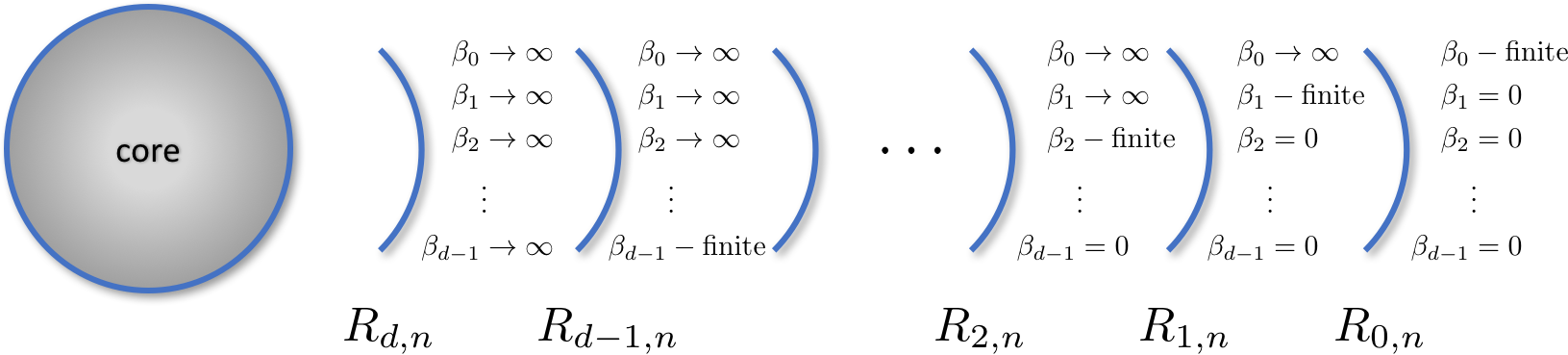}
\end{center}
\caption{\label{fig:crackle} The annuli described in \cite{adler2013crackle,owada_limit_2015}. Different homology degrees show up at different radii, where the lower degrees reach further away from the origin. Close to the origin we have a region called the `core' where the \cech complex is contractible.}
\end{figure}

The general setup in \cite{adler2013crackle,owada_limit_2015}  is the following. Let $f:\R^d\to\R$ be a probability distribution function whose support is $\R^d$, and let $\cC_r(n)$ defined as before. The results in these paper show that as $n\to\infty$ and $r\to 0$, even when $\Lambda \gg \log n$, many cycles can still show up far away from the origin. Moreover, it can be shown that homology has a very organized spatial structure.
Loosely speaking, we can split $\R^d$ into a sequence of annuli, such that inside each annulus we can find connected components that generate homology at different degrees.
More concretely - there is a sequence of radii $R_{0,n} > R_{1,n} > R_{2,n} \cdots R_{d,n}$ (depending on $r$ and $f$) such that inside the annulus $(R_{k,n}, R_{k-1,n})$ we have that $\beta_k$ is finite, $\beta_i \to\infty$ for $i < k$ and $\beta_i \to 0$ for $i > k$ (where by $\beta_i$ we mean the dimension of the $i$-th homology generated on connected components made of vertices that are contained in the specified annulus). In addition, there is a smaller radius $R_{c,n} < R_{d,n}$ such that the \cech complex inside $B_{R_{c,n}}(0)$ is contractible, and thus contains no nontrivial homology. This region is referred to as `the core'. This phenomenon is described in Figure \ref{fig:crackle}.

The work in \cite{owada_limit_2015} studies this phenomena in detail, discussing the differences between light and heavy-tailed distributions, and proving that there is a limiting Poisson law that describes the spatial distribution of cycles  appearing in each annulus.

\section{Persistent homology} \label{sec:persistent}

Persistent homology is one of the most heavily used tools in applied topology, or TDA (cf. \cite{Carlsson09,Ghrist08}). However, very little is known about its probabilistic properties. Briefly, the persistent homology of a \cech or a Rips complex tracks the evolution of the homology of the complex as the radius $r$ changes from zero to infinity. 
In this section we will review some recent work related to the persistent homology of random geometric complexes \cite{BKS14,duy_limit_2016}. 

\begin{figure}
\centering
\includegraphics[scale=0.5]{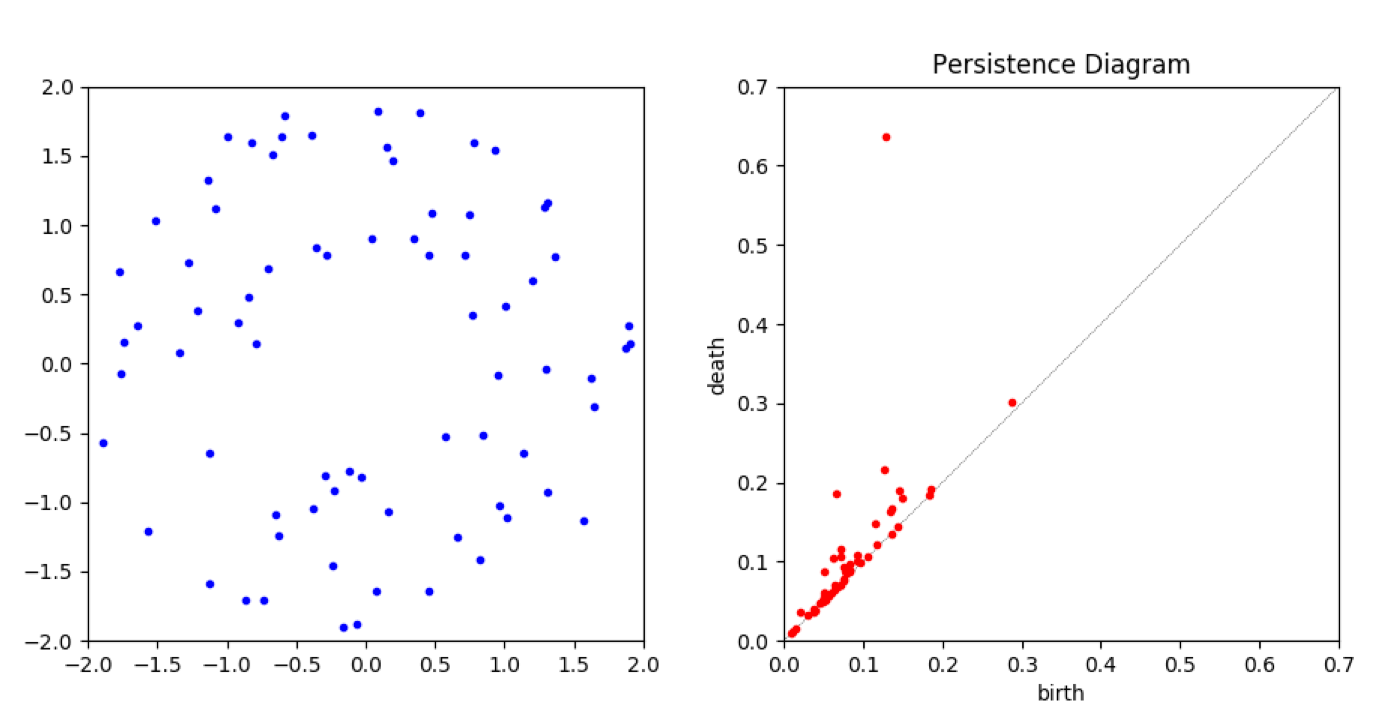}
\caption{\label{fig:pd} The persistence diagram of a random \cech filtration. The point process (on the left) is generated on an annulus in $\R^2$. The $H_1$ persistence diagram (on the right) describes the birth and death times (radii) of all the $1$-cycles that appear in this filtration. Notice that most of the points in the persistence diagram are close to the diagonal (where death=birth), and one might consider these cycles as "noise".
There is one point that stands out in the diagram, which corresponds to the hole of the annulus. The persistent homology was computed using the GUDHI library \cite{gudhi:urm}.
}
\end{figure}

Loosely speaking, the $k$-th persistent homology $\PH_k$ contains a list of all the $k$-dimensional nontrivial cycles that are created (and later terminated) in a geometric complex as $r$ is increased from $0$ to $\infty$. For every cycle $\gamma\in \PH_k$, we can assign a pair of values $(\gamma_{\birth}, \gamma_{\death})$ that represent the radii at which $\gamma$ appear and vanish (born and dies), respectively. A popular way to visualize the information provided by persistent homology is called the \emph{persistence diagram}.
Here, for every cycle $\gamma \in \PH_k$ we place a single point in the plane, where the $x$ and $y$ axes correspond to the birth and death times, respectively. Figure \ref{fig:pd} shows the persistence diagram of $H_1$ for a random \cech filtration.

\subsection{Limit theorems for persistence diagrams}
Denote by $\xi_k$ the persistence diagram for $\PH_k$. Clearly, $\xi_k \subset \Delta$, where $\Delta := \{ (x,y) : 0 \le x < y < \infty\}$, since death always occurs after birth (see Figure \ref{fig:pd}).

In \cite{duy_limit_2016}, the \cech and the Rips complex were considered, taken over stationary point processes $\Phi$
(as discussed in Section \ref{sec:stationary}). In this case, taking $\xi_{k,n}$ to be the $k$-th persistence diagram of $\Phi_n$, then $\xi_{k,n}$ is a random point process, or random Radon measure,  in $\R^2$. One of the main theorems in \cite{duy_limit_2016} states that as $n\to\infty$ this measure has a nonrandom limit $\nu_k$. In particular,
\begin{theorem}[Dul et. al., \cite{duy_limit_2016}]
If $\Phi$ is a stationary point process in $\R^d$ with finite moments, then there exists a unique Radon measure $\nu_k$ on $\Delta$ such that
\[
	\frac{1}{n^d}\mean{\xi_{k,n}} \xrightarrow{n\to\infty} \nu_k,
\]
where the convergence is in terms of the vague convergence of measures on $\Delta$. If, in addition, $\Phi$ is ergodic, then almost surely
\[
	\frac{1}{n^d}{\xi_{k,n}} \xrightarrow{n\to\infty} \nu_k,
\]
\end{theorem}
Under some additional conditions on $\Phi$ they show that the support of the limiting measure $\nu_k$ is the subspace $R_k\subset \Delta$ of all (birth, death) pairs realizable by the corresponding filtration (which can be \cech, Rips, and others). For example, for the \cech filtration
\[
	R_k = \begin{cases} \set{0}\times (0,\infty] & k = 0,\\
	\Delta & 1\le k \le d-1,\\
	\emptyset & k \ge d.
	\end{cases}
\]
In addition to the convergence of the entire measure, they study the variables $\beta_k^{r,s}$ counting cycles with $\gamma_{\birth} \le r$ and $\gamma_{\death} \ge s$.  Using similar techniques to the ones in \cite{YSA14} they prove a law of large numbers and a central limit theorem.

\subsection{Maximal cycles in persistent homology}
In this section we review the result in \cite{BKS14}, related to extremal cycles.
Traditionally, the \emph{persistence} (or significance) of a cycle $\gamma$ is measured by the difference $\gamma_{death}-\gamma_{birth}$.
In this work, persistence was measured by the ratio $\pi(\gamma) := \gamma_{death}/\gamma_{birth}$. There are a number of reasons to measure the persistence  of a cycle multiplicatively.

\begin{itemize}
\item The persistence measured this way is scale invariant, i.e.\ the persistence of cycles for $n$ points chosen uniformly in a cube $[0,1]^d$ will have the same distribution as for $n$ points chosen uniformly in a cube $[0, \lambda]^d$ for any $\lambda > 0$.
\item In a random geometric setting, one issue with measuring persistence by $\gamma_{death}-\gamma_{birth}$ is that both terms are tending to zero as the number of vertices goes to infinity, and $\gamma_{birth} \ll \gamma_{death}$. For the prominent cycles, $\gamma_{birth}\to 0$ much faster than $\gamma_{death}$, and therefore if we measure persistence as $\gamma_{death} - \gamma_{birth}$, then $\gamma_{birth}$ will just be a small error term and it will be hard to differentiate between them. The multiplicative way of measuring persistence is more informative.
\item Both \cech complexes $\C_r(n)$ and Vietoris--Rips  complexes $\cR_r(n)$ are central to the theory of persistent homology, and it is important to be able to compare them. The standard way of relating them is via the inclusion maps
$$\dots \hookrightarrow \C_r(n) \hookrightarrow \cR_r(n) \hookrightarrow \C_{\sqrt{2}r}(n) \hookrightarrow \cR_{\sqrt{2}r}(n) \hookrightarrow \dots$$
(In general $\C_r(n) \hookrightarrow \cR_r(n) \hookrightarrow \C_{\alpha r}(n)$ for \cech and Vietoris--Rips complexes in Euclidean space $\R^d$, as long as $\alpha \ge \sqrt{ 2d / (d+1)}$, as shown in Theorem 2.5 of \cite{dSG07}.)

So one may relate persistent homology between the two types of complexes. Because this relationship is naturally multiplicative in $r$, our results are stated in a way that holds for both types of complexes.
\end{itemize}
The result in \cite{BKS14}  was proven for a homogeneous Poisson process on the unit cube $[0,1]^d$. However, similar results should hold for any measurable density function $f$ on any $d$-dimensional compact and convex body, provided that $f$ is bounded from below and above. 

\begin{theorem} 
Let $\cP_n$ be a unit-intensity Poisson process on the unit cube $[0,1]^d$. Let $\PH_k(n)$ be the $k$-th dimensional persistent homology of either the \cech or the Rips filtration generated by $\cP_n$.
Define,
\[
	\Pi_k(n) := \max_{\gamma\in \PH_k(n)} \pi(\gamma),
\]
i.e.~$\Pi_k(n)$ is the maximal persistence of all $k$-cycles.
Then a.a.s.~we have that
\[
	\Pi_k(n) = \Theta\param{\param{\frac{\log n}{\log \log n} }^{1/k}}.
\]
The implied constants in the asymptotic notation $\Theta$ only depend on the underlying probability distribution.
\end{theorem}

Persistent homology is becoming a very popular and powerful data analysis tool. Studying this type of extremal behavior for persistent homology can be later used to provide a statistical analysis to persistent homology. For example, suppose that the data are sampled from a distribution supported on a manifold $M$ with non trivial homology that we wish to recover. Knowing the distribution of $\Pi_k$ for convex bodies (where homology is trivial), would enable us to develop statistical tests to differentiate between the signal (real cycles of $M$) and noise (artifacts of the sampling mechanism) in this type of data analysis problem.
Persistent homology in random contexts was studied earlier by Bubenik and Kim in \cite{BK07}.

\section{Open problems / future directions}\label{sec:future}

We close by mentioning several possible directions for future research.

\begin{itemize}
\item {\bf Sharper results in the thermodynamic limit.} Proving strong results for expectation of Betti numbers in the critical regime remains a challenging problem. The best result so far is that
$$\frac{\expect[\beta_k(n)] }{n} \to C,$$
where $C > 0$ is some constant which depends on the underlying distribution on $\R^d$ and the degree $k$ \cite{YSA14}. It would be a breakthrough to write an explicit formula for $C$ and we expect that the results would find applications in TDA.
\item {\bf Connections between the various models.} Is there a model for random geometric complex which approximates the sub-level sets of the Gaussian random field? See \cite{AT07} and \cite{AT11} for introduction and overview of Gaussian random fields and their topological properties. 
\item {\bf Torsion.} All of the results in this survey for homology of random geometric complexes do not depend on the choice of coefficients. In dimensions $d \ge 4$ and higher, these complexes will likely have torsion in integer homology. What can be said about the limiting distribution of this torsion group?
\item {\bf Higher-dimensional percolation theory.} All of the random geometric complexes discussed here are analogues of random geometric graphs where the number of vertices $n$ is finite and $n \to \infty$. Percolation theory is of a somewhat different flavor---one considers an infinite random graph, by taking a random subgraph of a lattice, and then analyzes large-scale structure such as whether or not an infinite connected component appears. Analogous lattice models with higher-dimensional cells have been studied, for example ``plaquette percolation.'' \cite{ACCFR83,GH10}. So rather than study homology-vanishing thresholds for finite random geometric complexes with size tending to infinity, one might study the appearance of ``infinite'' cycles in lattice models. So far, this seems to be relatively unexplored.
\end{itemize}

\vspace{20pt}
\noindent On behalf of all authors, the corresponding author states that there is no conflict of interest. 
\vspace{10pt}

\bibliography{survey}


\end{document}